\documentclass{article}
\usepackage{amssymb}
\usepackage{amsmath}
\usepackage{graphicx}

\input{tcilatex}
\begin{document}

\begin{center}
\ \ \ \ \ \ \ \ \ \ \ \ \ \ \ \ \ \ \ \ \ \ \ \ \ \ \ \ \ \ \ \ \ \ \ \ \ \
\ \ \ \ \ \ \ \ \ \ \ \ \ 

\bigskip\ \ \ \textbf{Veli B. Shakhmurov}\ \ \ 

Okan University, Department of Mechanical Engineering, Akfirat, Tuzla 34959
Istanbul, Turkey, E-mail: veli.sahmurov@okan.edu.tr

\bigskip \textbf{Separability properties of singular degenerate abstract
differential operators and applications\ }

\ \ \ \ \ \ \ \ \ \ \ \ \ \ \ \ \ \ \ \ 

\ \ \ \ \ \ \ \ \ \ \ \ \ \ \ \ \ \ \ \ \ \ \ 

\textbf{AMS: 34G10, 35J25, 35J70\ \ }\ \ \ \ \ \ \ \ \ \ \ \ \ 

\textbf{Abstract}
\end{center}

In this paper, we study the separability and spectral properties of singular
degenerate elliptic equations in vector valued $L_{p}$ spaces. We prove that
a realization operator by this equation with some boundary conditions is
separable and Fredholm in $L_{p}$. The leading part of the associated
differential operator is not self-adjoint. The sharpe estimate of the
resolvent, discreetness of spectrum and completeness of root elements of
this operator is obtained. Moreover, we show that this operator is positive
and generates a holomorphic $C_{0}$-semigroups on $L_{p}.$ In application,
we examine the regularity properties of degenerate elliptic problem with
Wentzell--Robin boundary conditions and boundary value problem for system of
degenerate elliptic equations of either finite or infinite number.


\begin{center}
\ \ \textbf{1. Introduction, notations and background }
\end{center}

In this work, boundary value problem (BVP) for singular degenerate abstract
elliptic equations are considered. BVPs for abstract differential equations
(ADEs) have been studied extensively by many researchers (see e.g. $\left[
1-3\right] ,$ $\left[ 6-8\right] $, $\left[ 10-19\right] $, $\left[ 21-22%
\right] $ and the references therein). A comprehensive introduction to the
ADEs and historical references may be found in $\left[ 1\right] $ and $\left[
22\right] .$ The maximal regularity properties for differential operator
equations have been investigated e.g. in $\left[ 6-8\right] $, $\left[ 14-18%
\right] $ and $\left[ 21-22\right] $. The main objective of the present
paper is to discuss the BVP for the following singular degenerate DOE 
\begin{equation}
\dsum\limits_{k=1}^{n}-x_{k}^{2\alpha _{k}}\frac{\partial ^{2}u}{\partial
x_{k}^{2}}+\dsum\limits_{k=1}^{n}x_{k}^{\alpha _{k}}A_{k}\left( x\right) 
\frac{\partial u}{\partial x_{k}}+Au=f\left( x\right) ,  \tag{1.1}
\end{equation}%
where $A$, $A_{k}$ are linear operators in a Banach space $E.$

We derive $L_{p}-$separability properties and sharp resolvent estimates of
the associated differential operator. Especially, we show that this
differential operator is $R$-positive and\ also is a generator of an
analytic semigroup.

By using separability properties of the elliptic problem $\left( 1.1\right) $
we derive spectral properties of differential operator $Q$ generated by $%
\left( 1.1\right) .$ Namely, we prove that the operator $Q$ is Fredholm in $%
L_{p}$, the inverse $Q^{-1}$ belong to some Schatten class $\sigma
_{q}\left( L_{p}\right) $ and the system of root functions of this operator
is complete in $L_{p}.$

\bigskip One of the most important aspects of this ADE considered here is
that the degeneration in different directions is at different speeds, in
general. Unlike the regular degenerate equations, because of the singularity
of the degeneracy of the equation, the boundary conditions are only given on
the lines\ without degeneracy.

In application, the BVP for infinity system of singular degenerate partial
differential equations and Wentzell-Robin type BVP for singular degenerate
partial differential equations on cylindrical domain are studied.

Since the Banach space $E$ is arbitrary and $A$ is a possible linear
operator, by choosing $E$ and $A$ we can obtain numerous classis of
degenerate elliptic and qusielliptic equations which have a different
applications. Let we choose $E=L_{2}\left( 0,1\right) $ and $A$ to be
differential operator providing the Wentzell-Robin boundary condition
defined by 
\begin{equation*}
D\left( A\right) =\left\{ u\in W_{2}^{2}\left( 0,1\right) ,\text{ }A\left(
j\right) u\left( j\right) =0,\text{ }j=0,1\right\} ,\text{ }
\end{equation*}%
\begin{equation*}
\text{ }Au=a\left( y\right) u^{\left( 2\right) }+b\left( y\right) u^{\left(
1\right) }\text{ for all }y\in \left( 0,1\right) ,
\end{equation*}%
where $a$ is positive and $b$ is a real-valued functions on $\left(
0,1\right) $. By virtue of $L_{p}-$ regularity properties of $\left(
1.1\right) $ (see Theorem 2.1)\ we obtain the separability properties of
Wentzell-Robin type BVP for singular degenerate elliptic equation

\begin{equation}
\dsum\limits_{k=1}^{n}x_{k}^{2\alpha _{k}}\frac{\partial ^{2}u}{\partial
x_{k}^{2}}+a\left( y\right) \frac{\partial ^{2}u}{\partial y^{2}}+b\left(
y\right) \frac{\partial u}{\partial y}=f\left( x,y\right) ,\text{ }x\in G%
\text{, }y\in \left( 0,1\right) ,  \tag{1.2}
\end{equation}%
\ \ \ 

\begin{equation*}
L_{k}u=0,\text{ for a.e. }y\in \left( 0,1\right) ,\text{ }
\end{equation*}%
\begin{equation}
A\left( j\right) u\left( j\right) =a\left( j\right) u_{yy}\left( x,j\right)
+b\left( j\right) u_{y}\left( x,j\right) =0\text{, }  \tag{1.3}
\end{equation}%
\begin{equation*}
j=0,1,\text{ for a.e. }x\in G,
\end{equation*}%
in the mixed $L_{\mathbf{p}}\left( \Omega \right) $ spaces,\ where $L_{k}$
are boundary conditions with respect $x\in G\subset R^{n}$ that will be
definet in late and $L_{\mathbf{p}}\left( \Omega \right) $ denotes the space
of all $\mathbf{p}$-summable complex-valued\ functions with mixed norm and \ 
\begin{equation*}
\Omega =G\times \left( 0,1\right) ,\text{ }\mathbf{p=}\left( p,2\right) .
\end{equation*}

Note that, the regularity properties of Wentzell-Robin type problems for
elliptic and parabolic equations were studied e.g. in $\left[ \text{5, 9, 11}%
\right] $ and the references therein.

Let $\gamma =\gamma \left( x\right) $ be a positive measurable function on a
domain $\Omega \subset R^{n}.$ Here, $L_{p,\gamma }\left( \Omega ;E\right) $
denote the space of strongly measurable $E$-valued functions that are
defined on $\Omega $ with the norm

\begin{equation*}
\left\Vert f\right\Vert _{L_{p,\gamma }}=\left\Vert f\right\Vert
_{L_{p,\gamma }\left( \Omega ;E\right) }=\left( \int \left\Vert f\left(
x\right) \right\Vert _{E}^{p}\gamma \left( x\right) dx\right) ^{\frac{1}{p}},%
\text{ }1\leq p<\infty .
\end{equation*}

For $\gamma \left( x\right) \equiv 1$ the space $L_{p,\gamma }\left( \Omega
;E\right) $ will be denoted by $L_{p}=L_{p}\left( \Omega ;E\right) .$

\ The Banach space\ $E$ is called an $UMD$-space if\ the Hilbert operator $%
\left( Hf\right) \left( x\right) =\lim\limits_{\varepsilon \rightarrow
0}\int\limits_{\left\vert x-y\right\vert >\varepsilon }\frac{f\left(
y\right) }{x-y}dy$ \ is bounded in $L_{p}\left( R,E\right) ,$ $p\in \left(
1,\infty \right) $ (see. e.g. $\left[ 4\right] $). $UMD$ spaces include e.g. 
$L_{p}$, $l_{p}$ spaces and Lorentz spaces $L_{pq},$ $p$, $q\in \left(
1,\infty \right) $.

Let $\mathbb{C}$ be the set of the complex numbers and\ 
\begin{equation*}
S_{\varphi }=\left\{ \lambda ;\text{ \ }\lambda \in \mathbb{C}\text{, }%
\left\vert \arg \lambda \right\vert \leq \varphi \right\} \cup \left\{
0\right\} ,\text{ }0\leq \varphi <\pi .
\end{equation*}

\ Let $E_{1}$ and $E_{2}$ be two Banach spaces. $L\left( E_{1},E_{2}\right) $
denotes the space of bounded linear operators from $E_{1}$ into $E_{2}.$ For 
$E_{1}=E_{2}=E$ it will be denoted by $L\left( E\right) .$

\textbf{Definition 1. }A linear operator\ $A$ is said to be $\varphi $%
-positive in a Banach\ space $E$ with bound $M>0$ if $D\left( A\right) $ is
dense on $E$ and $\left\Vert \left( A+\lambda I\right) ^{-1}\right\Vert
_{L\left( E\right) }\leq M\left( 1+\left\vert \lambda \right\vert \right)
^{-1}$ $\ $for any $\lambda \in S_{\varphi },$ $0\leq \varphi <\pi ,$ where $%
I$ is an identity operator in $E$. Sometimes $A+\lambda I$\ will be denoted
by $A+\lambda $ or$\ A_{\lambda }$. It is known $\left[ \text{20, \S 1.15.1}%
\right] $ that a positive operator $A$ has well-defined fractional powers\ $%
A^{\theta }.$

\textbf{Remark 1.1. }By virtue of $\left[ \text{20, \S\ 1.13}\right] $ if $A$
is $\varphi $-positive in $E$, then the operator $-A^{\alpha }$ generate an
analytic semigroup $U_{A^{\alpha }}\left( t\right) $ for $0<\alpha \leq 1$, $%
\varphi \geq \frac{\pi }{2}$ and for $\alpha \leq \frac{1}{2}$, $\varphi <%
\frac{\pi }{2}.$ Moreover, there exists a positive constant $\omega $ such
that the estimate holds 
\begin{equation*}
\left\Vert U_{A^{\alpha }}\left( t\right) \right\Vert _{L\left( E\right)
}\leq Me^{-\omega t}.
\end{equation*}

Let $E\left( A^{\theta }\right) $ denote the space $D\left( A^{\theta
}\right) $ equipped with the norm 
\begin{equation*}
\left\Vert u\right\Vert _{E\left( A^{\theta }\right) }=\left( \left\Vert
u\right\Vert ^{p}+\left\Vert A^{\theta }u\right\Vert ^{p}\right) ^{\frac{1}{p%
}},1\leq p<\infty ,\text{ }0<\theta <\infty .
\end{equation*}

Let $E_{1}$ and $E_{2}$ be two Banach spaces. Now $\left( E_{1},E_{2}\right)
_{\theta ,p}$, $0<\theta <1,1\leq p\leq \infty $ will denote interpolation
spaces obtained from $\left\{ E_{1},E_{2}\right\} $ by the $K$ method \ $%
\left[ \text{20, \S 1.3.1}\right] $.

\textbf{Definition 2. }Let $\mathbb{N}$ denote the set of natural numbers
and $\left\{ r_{j}\right\} $ is a sequence of independent symmetric $\left\{
-1,1\right\} $-valued random variables on $\left[ 0,1\right] $. A set $%
K\subset L\left( E_{1},E_{2}\right) $ is called $R$-bounded if there is a
positive constant $C$ such that for all $T_{1},T_{2},...,T_{m}\in K$ and $%
u_{1,}u_{2},...,u_{m}\in E_{1},$ $m\in \mathbb{N}$ 
\begin{equation*}
\int\limits_{0}^{1}\left\Vert \sum\limits_{j=1}^{m}r_{j}\left( y\right)
T_{j}u_{j}\right\Vert _{E_{2}}dy\leq C\int\limits_{0}^{1}\left\Vert
\sum\limits_{j=1}^{m}r_{j}\left( y\right) u_{j}\right\Vert _{E_{1}}dy.
\end{equation*}

The smallest $C$ for which the above estimate holds is called a $R$-bound of
the collection $K$ and denoted by $R\left( K\right) .$

\textbf{Definition 3.} The $\varphi $-positive operator $A$ is said to be $R$%
-positive in $E$ if the set $\left\{ \lambda \left( A+\lambda I\right) ^{-1}%
\text{: }\lambda \in S_{\varphi }\right\} ,$ $0\leq \varphi <\pi $ is $R$%
-bounded.

Let $E_{1}$ and $E_{2}$ be two Banach spaces. $\sigma _{\infty }\left(
E_{1},E_{2}\right) $ denotes the space of all compact operators from $E_{1}$
to $E_{2}.$ For $E_{1}=E_{2}=E$ it will be denoted by $\sigma _{\infty
}\left( E\right) .$

$s_{j}\left( A\right) $ will denote approximation numbers of operator $A$ $%
\left[ \text{20, \S\ 1.16.1}\right] $. Let 
\begin{equation*}
\sigma _{q}\left( E_{1},E_{2}\right) =\left\{ A\in \sigma _{\infty }\left(
E_{1},E_{2}\right) ,\text{ }\sum\limits_{j=1}^{\infty }s_{j}^{q}\left(
A\right) <\infty ,\text{ }1\leq q<\infty \right\} .
\end{equation*}

Here, $\Omega $ is a domain in $R^{n}$. Assume $E_{0}$ and $E$ are two
Banach spaces so that $E_{0}$ is continuously and densely embedded into $E$.
Let $\gamma _{k}=\gamma _{k}\left( x\right) $ be a positive measurable
functions on $\Omega $ and $\gamma =\left( \gamma _{1},\gamma
_{2},...,\gamma _{n}\right) $. Consider, the Sobolev-Lions type space\ $%
W_{p,\gamma }^{m}\left( \Omega ;E_{0},E\right) $, i.e. the space consisting
of all functions $u\in L_{p}\left( \Omega ;E_{0}\right) $ that have
generalized derivatives $D_{k}^{m}u=\frac{\partial ^{m}u}{\partial x_{k}^{m}}%
\in L_{p,\gamma _{k}}\left( \Omega ;E\right) $ equipped with the norm 
\begin{equation*}
\ \left\Vert u\right\Vert _{W_{p,\gamma }^{m}\left( \Omega ;E_{0},E\right)
}=\left\Vert u\right\Vert _{L_{p}\left( \Omega ;E_{0}\right)
}+\dsum\limits_{k=1}^{n}\left\Vert \gamma _{k}^{m}\frac{\partial ^{m}u}{%
\partial x_{k}^{m}}\right\Vert _{L_{p}\left( \Omega ;E\right) }<\infty .
\end{equation*}

Let $\chi =\chi \left( t\right) $ be a positive measurable function on $%
\left( 0,a\right) $ and 
\begin{equation*}
u^{\left[ i\right] }\left( t\right) =\left( \chi \left( t\right) \frac{d}{dt}%
\right) ^{i}u\left( t\right) .
\end{equation*}%
Consider the following $E-$valued weighted function spaces 
\begin{equation*}
W_{p,\chi }^{\left[ m\right] }\left( 0,a;E_{0},E\right) =\left\{ u;u\in
L_{p}\left( 0,a;E_{0}\right) \right. ,\ u^{\left[ m\right] }\in L_{p}\left(
0,a;E\right) ,
\end{equation*}

\begin{equation*}
\left\Vert u\right\Vert _{W_{p,\chi }^{\left[ m\right] }}=\left. \left\Vert
u\right\Vert _{L_{p}\left( 0,a;E_{0}\right) }+\left\Vert u^{\left[ m\right]
}\right\Vert _{L_{p}\left( 0,a;E\right) }<\infty \right\} .
\end{equation*}

\begin{equation*}
W_{p,\chi }^{m}\left( 0,a;E_{0},E\right) =\left\{ u;u\in L_{p,\chi }\left(
0,a;E_{0}\right) \right. ,\ u^{\left( m\right) }\in L_{p,\chi }\left(
0,a;E\right) ,
\end{equation*}

\begin{equation*}
\left\Vert u\right\Vert _{W_{p,\chi }^{m}}=\left. \left\Vert u\right\Vert
_{L_{p,\chi }\left( 0,a;E_{0}\right) }+\left\Vert u^{\left( m\right)
}\right\Vert _{L_{p,\chi }\left( 0,a;E\right) }<\infty \right\} .
\end{equation*}

Let 
\begin{equation*}
\alpha =\left( \alpha _{1},\alpha _{2},...,\alpha _{n}\right) \text{, }D^{%
\left[ \alpha \right] }=D_{1}^{\left[ \alpha _{1}\right] }D_{2}^{\left[
\alpha _{2}\right] }...D_{n}^{\left[ \alpha _{n}\right] },\text{ }D_{k}^{%
\left[ i\right] }.=\left( \gamma _{k}\left( x\right) \frac{\partial }{%
\partial x_{k}}\right) ^{i}.
\end{equation*}

Consider the space\ $W_{p,\gamma }^{\left[ m\right] }\left( \Omega
;E_{0},E\right) $, consisting of all functions $u\in L_{p}\left( \Omega
;E_{0}\right) $ that have generalized derivatives $D_{k}^{\left[ m\right]
}u\in L_{p,}\left( \Omega ;E\right) $ with the norm 
\begin{equation*}
\ \left\Vert u\right\Vert _{W_{p,\gamma }^{\left[ m\right] }\left( \Omega
;E_{0},E\right) }=\left\Vert u\right\Vert _{L_{p}\left( \Omega ;E_{0}\right)
}+\dsum\limits_{k=1}^{n}\left\Vert D_{k}^{\left[ m\right] }u\right\Vert
_{L_{p}\left( \Omega ;E\right) }<\infty .
\end{equation*}

From $\left[ \text{15, Theorem 1, Theorem 3}\right] $ we obtain

\textbf{Theorem A}$_{1}$. Suppose the following conditions are satisfied:

(1) $E$ is an UMD space and\ $A$ is an $R$-positive operator in $E;$

(3) $\gamma _{k}\left( x\right) =x_{k}^{\nu _{k}}$, $\nu _{k}\in \left(
1,p\right) $, $p\in \left( 1,\infty \right) $, $m$ is an integer and $0\leq
\mu \leq 1-\frac{\left\vert \alpha \right\vert }{m},$ $1<p<\infty $;

(4)\ $\Omega =\dprod\limits_{k=1}^{n}\left( 0,a_{k}\right) .$

Then, the embedding 
\begin{equation*}
D^{\alpha }W_{p,\gamma }^{m}\left( \Omega ;E\left( A\right) ,E\right)
\subset L_{p}\left( \Omega ;E\left( A^{1-\frac{\left\vert \alpha \right\vert 
}{m}-\mu }\right) \right)
\end{equation*}
is continuous. Moreover for all $h>0$ with $h\leq $ $h_{0}<\infty $ and $%
u\in W_{p,\gamma }^{m}\left( \Omega ;E\left( A\right) ,E\right) $ the
following uniform estimate holds 
\begin{equation*}
\left\Vert D^{\alpha }u\right\Vert _{L_{p,\gamma _{\alpha }}\left( \Omega
;E\left( A^{1-\frac{\left\vert \alpha \right\vert }{m}-\mu }\right) \right)
}\leq h^{\mu }\left\Vert u\right\Vert _{W_{p,\gamma }^{m}\left( \Omega
;E\left( A\right) ,E\right) }+h^{-\left( 1-\mu \right) }\left\Vert
u\right\Vert _{L_{p}\left( \Omega ;E\right) },
\end{equation*}%
where 
\begin{equation*}
\gamma _{\alpha }\left( x\right) =\dprod\limits_{k=1}^{n}x_{k}^{\alpha
_{k}\nu _{k}}.
\end{equation*}

\textbf{Theorem A}$_{2}$. Assume the conditions of Theorem A$_{1}$ are
satisfied. Moreover, suppose $a_{k}<\infty $ and $A^{-1}$ is a compact
operator $E.$ Then for $0<\mu \leq 1-\frac{\left\vert \alpha \right\vert }{m}
$ the embedding 
\begin{equation*}
D^{\alpha }W_{p,\gamma }^{m}\left( \Omega ;E\left( A\right) ,E\right)
\subset L_{p}\left( \Omega ;E\left( A^{1-\varkappa -\mu }\right) \right)
\end{equation*}
is compact.

Let $I\left( E_{0},E\right) $ denote the embedding operator from $E_{0}$ to $%
E.$ By reasoning as in $\left[ \text{14, Theorem 3.1}\right] $ we have

\textbf{Theorem A}$_{3}.$ Let $E$ be Banach spaces with base $\alpha _{k}\in
\left( 1,p\right) $ for $p\in \left( 1,\infty \right) $ and $\alpha _{k}<m$.
Suppose the embedding $E_{0}\subset E$ is compact and 
\begin{equation*}
s_{j}\left( I\left( E_{0},E\right) \right) \sim j^{-\frac{1}{\nu }},\text{ }%
\nu >0,\text{ }j=1,2,...,\infty .
\end{equation*}%
Then 
\begin{equation*}
s_{j}\left( I\left( W_{p,\alpha }^{m}\left( G;E_{0},E\right) ,L_{p}\left(
G;E\right) \right) \right) \sim j^{-\frac{1}{\nu +\varkappa }}\text{, }%
\varkappa =\dsum\limits_{k=1}^{n}\frac{1}{m-\alpha _{k}}.
\end{equation*}

Consider the BVP 
\begin{equation}
-u^{\left( 2\right) }\left( t\right) +Au\left( t\right) =f\left( t\right) , 
\tag{1.4}
\end{equation}%
\begin{equation*}
L_{1}u=\sum\limits_{i=0}^{m}\left[ \delta _{i}u^{\left( i\right) }\left(
a\right) +\sum\limits_{j=1}^{N}\nu _{ij}u^{\left( i\right) }\left(
t_{ij}\right) \right] =0,
\end{equation*}%
where $m\in \left\{ 0,1\right\} ;$ $\delta _{i},$ $\nu _{ij},$ are complex
numbers, $t_{ij}\in \left( 0,a\right) $ and $A$ is a linear operator in $E.$

\textbf{Cond\i t\i on 1.1. }Let the following conditions be satisfied:

(1)\ $E$ is a UMD\ space and $A$ is a $R$ positive operator in $E$;

(2) $\delta _{m}\neq 0$, $\gamma \left( t\right) =t^{\nu }$, $\nu \in \left(
1,p\right) $ for $p\in \left( 1,\infty \right) $;

(3) Here, $t_{0}=\min\limits_{j}t_{1j}$, $U_{A^{\frac{1}{2}}}^{-1}\left(
t_{0}\right) \in L\left( E\right) $ and 
\begin{equation*}
\sum\limits_{i=0}^{m}\left[ \left\vert \delta _{i}\right\vert e^{-\omega
\left( a-t_{0}\right) }+\sum\limits_{j=1}^{N}\left\vert \nu _{ij}\right\vert
e^{-\omega \left( t_{ij}-t_{0}\right) }\right] <\left\vert \nu
_{0}\right\vert ,
\end{equation*}%
where $\omega $ is a positive constant defined in the Remark 1.1.

Let $\gamma \left( t\right) =t^{\frac{\nu }{1-\nu }}$. In a similar way as
in $\left[ \text{16, Theorem 5.1}\right] $ we obtain

\textbf{Theorem A}$_{4}$. Assume the Condition 1.1 are satisfied. Then, the
problem $\left( 1.4\right) $ has a unique solution%
\begin{equation*}
u\in W_{p,\gamma }^{2}\left( a,\infty ;E\left( A\right) ,E\right)
\end{equation*}
for all \ $f\in L_{p,\gamma }\left( a,\infty ;E\right) $, $\left\vert \arg
\lambda \right\vert \leq \varphi $ with sufficiently large $\left\vert
\lambda \right\vert $ and the uniform coercive estimate holds

\begin{equation*}
\sum\limits_{i=0}^{2}\left\vert \lambda \right\vert ^{1-\frac{i}{2}%
}\left\Vert u^{\left( i\right) }\right\Vert _{L_{p,\gamma }\left( a,\infty
;E\right) }+\left\Vert Au\right\Vert _{L_{p,\gamma }\left( a,\infty
;E\right) }\leq C\left\Vert f\right\Vert _{L_{p,\gamma }\left( a,\infty
;E\right) }.
\end{equation*}

\begin{center}
\bigskip \textbf{2. Singular degenerate abstract elliptic equations}
\end{center}

\ Consider the BVP for the following singular degenerate ADO 
\begin{equation}
\dsum\limits_{k=1}^{n}\left[ -x_{k}^{2\alpha _{k}}\frac{\partial ^{2}u}{%
\partial x_{k}^{2}}+x_{k}^{\alpha _{k}}A_{k}\left( x\right) \frac{\partial u%
}{\partial x_{k}}\right] +Au+\lambda u=f\left( x\right) ,\text{ }x\in G, 
\tag{2.1}
\end{equation}

\begin{equation}
L_{k}u=\sum\limits_{i=0}^{m_{k}}\left[ \delta _{ki}u_{x_{k}}^{\left[ i\right]
}\left( a_{k},x\left( k\right) \right) +\sum\limits_{j=0}^{N_{k}}\nu
_{kij}u_{x_{k}}^{\left[ i\right] }\left( x_{kij},x\left( k\right) \right) %
\right] =0\text{, }  \tag{2.2}
\end{equation}%
where $x\left( k\right) \in G_{k}$, $x_{kij}\in \left( 0,a_{k}\right) $ and%
\begin{eqnarray*}
u_{x_{k}}^{\left[ i\right] } &=&\left[ x_{k}^{\alpha _{k}}\frac{\partial }{%
\partial x_{k}}\right] ^{i}u\left( x\right) ,\text{ }G=\dprod%
\limits_{k=1}^{n}\left( 0,a_{k}\right) ,\text{ }G_{k}=\dprod\limits_{j\neq
k}\left( 0,a_{j}\right) , \\
m_{k} &\in &\left\{ 0,1\right\} \text{, }x\left( k\right) =\left(
x_{1},x_{2},...,x_{k-1},x_{k+1},...,x_{n}\right) \text{, }j,\text{ }%
k=1,2,...,n;
\end{eqnarray*}%
$\delta _{ki},$ $\nu _{kij}$ are complex numbers, $\lambda $ is a complex
parameter, $A$ and $A_{k}\left( x\right) $ are linear operators in a Banach
space $E.$

Let we denote $W_{p,\gamma }^{m}\left( \Omega ;E\left( A\right) ,E\right) $
by $W_{p,\alpha }^{m}\left( \Omega ;E\left( A\right) ,E\right) $ for $\gamma
_{k}\left( x\right) =x_{k}^{\alpha }$ $.$

\textbf{Cond\i t\i on 2.1. }Assume the following conditions are satisfied:

(1) $E$ is an UMD space and $A$ is a $R$-positive operator in $E;$

(2) $\delta _{km_{k}}\neq 0$, $\alpha _{k}\in \left( 1,p\right) $ for $p\in
\left( 1,\infty \right) $ and $k=1,2,...,n;$

(3) Here, $x_{k0}=\min\limits_{j}x_{k1j}$ and $U_{A^{\frac{1}{2}%
}}^{-1}\left( x_{k0}\right) \in L\left( E\right) .$ Moreover, 
\begin{equation*}
\sum\limits_{i=0}^{m}\left[ \left\vert \delta _{ki}\right\vert e^{-\omega
\left( a_{k}-x_{k0}\right) }+\sum\limits_{j=1}^{N}\left\vert \nu
_{kij}\right\vert e^{-\omega \left( x_{kij}-x_{k0}\right) }\right]
<\left\vert \nu _{k0}\right\vert ,
\end{equation*}

where $\omega $ is a positive constant defined in the Remark 1.1.

Let $\alpha =\left( \alpha _{1},\alpha _{2},...,\alpha _{n}\right) ,$ $%
\gamma _{k}\left( x\right) =x_{k}^{\alpha _{k}}.$ The main result is the
following:

\textbf{Theorem 2.1. }Assume the Condition 2.1 are hold and for any $%
\varepsilon >0$ there is a positive constant $C\left( \varepsilon \right) $
such that%
\begin{equation*}
\left\Vert A_{k}\left( x\right) u\right\Vert \leq \varepsilon \left\Vert
u\right\Vert _{\left( E\left( A\right) ,E\right) _{\frac{1}{2},\infty
}}+C\left( \varepsilon \right) \left\Vert u\right\Vert \text{ for\ }u\in
\left( E\left( A\right) ,E\right) _{\frac{1}{2},\infty }.
\end{equation*}

Then, problem $\left( 2.1\right) -\left( 2.2\right) $ has a unique solution $%
u\in W_{p,\alpha }^{2}\left( G;E\left( A\right) ,E\right) $ for $f\in
L_{p}\left( G;E\right) $ and sufficiently large $\left\vert \lambda
\right\vert $ with $\left\vert \arg \lambda \right\vert \leq \varphi $ and
the following uniform coercive estimate holds

\begin{equation}
\sum\limits_{k=1}^{n}\sum\limits_{i=0}^{2}\left\vert \lambda \right\vert ^{1-%
\frac{i}{2}}\left\Vert x_{k}^{i\alpha _{k}}\frac{\partial ^{i}u}{\partial
x_{k}^{i}}\right\Vert _{L_{p}\left( G;E\right) }+\left\Vert Au\right\Vert
_{L_{p}\left( G;E\right) }\leq M\left\Vert f\right\Vert _{L_{p}\left(
G;E\right) }.  \tag{2.3}
\end{equation}

For proving the main theorem, consider at first the BVP for the singular
degenerate ordinary DOE

\begin{equation}
\ -u^{\left[ 2\right] }\left( t\right) +\left( A+\lambda \right) u\left(
t\right) =f,\text{ }t\in \left( 0,a\right) ,  \tag{2.4}
\end{equation}%
\begin{equation*}
L_{1}u=\sum\limits_{i=0}^{m}\left[ \delta _{i}u^{\left[ i\right] }\left(
a\right) +\sum\limits_{j=1}^{N}\nu _{ij}u^{\left[ i\right] }\left(
t_{ij}\right) \right] =0,
\end{equation*}%
where $u^{\left[ i\right] }=\left( t^{\nu }\frac{d}{dt}\right) ^{i}$, $\nu
>1,$ $m\in \left\{ 0,1\right\} ,$ $\delta _{i},$ $\nu _{ij},$ are complex
numbers, $t_{ij}\in \left( 0,a\right) $ and $A$ is a linear operator in $E.$

Let 
\begin{equation*}
\gamma =\gamma \left( t\right) =t^{\nu }.
\end{equation*}

\textbf{Remark 2.1. }Consider the following substitution 
\begin{equation}
\tau =-\int\limits_{t}^{a}\gamma ^{-1}\left( z\right) dz,\text{ }t=t\left(
\tau \right) =\left[ a^{1-\nu }-\left( \nu -1\right) \tau \right] ^{\frac{1}{%
1-\nu }},  \tag{2.5}
\end{equation}

Under the substitution $\left( 2.5\right) $ the spaces $L_{p}\left(
0,a;E\right) $, $W_{p,\gamma }^{\left[ 2\right] }\left( 0,a;E\left( A\right)
,E\right) $ are mapped isomorphically onto weighted spaces%
\begin{equation*}
L_{p,\tilde{\gamma}}\left( -\infty ,0;E\right) ,\text{ }W_{p,\tilde{\gamma}%
}^{2}\left( -\infty ,0;E\left( A\right) ,E\right) ,
\end{equation*}%
respectively, where 
\begin{equation*}
\tilde{\gamma}=\gamma \left( t\left( \tau \right) \right) =\left[ a^{1-\nu
}-\left( \nu -1\right) \tau \right] ^{\frac{\nu }{1-\nu }}.
\end{equation*}%
Moreover, under the substitution $\left( 2.5\right) $ the problem $\left(
2.1\right) -\left( 2.2\right) $ is transformed into the following
undegenerate problem

\begin{equation}
-u^{\left( 2\right) }\left( \tau \right) +Au\left( \tau \right) =\tilde{f}%
\left( \tau \right) ,  \tag{2.6}
\end{equation}%
\begin{equation*}
L_{1}u=\sum\limits_{i=0}^{m}\left[ \delta _{i}u^{\left( i\right) }\left(
0\right) +\sum\limits_{j=1}^{N}\nu _{ij}u^{\left( i\right) }\left( \tau
_{ij}\right) \right] =0
\end{equation*}%
considered in the weighted space $L_{p,\tilde{\gamma}}\left( -\infty
;0;E\right) ,$ where $\tilde{f}\left( \tau \right) =f\left( t\left( \tau
\right) \right) $ and $\tau _{ij}\in \left( -\infty ,0\right) .$

By using Theorem A$_{4}$ we have

\textbf{Proposition 2.}1\textbf{. }Assume the Condition 1.1 are satisfied
with $t_{1j}=\tau _{1j}$. Then, the problem $\left( 2.4\right) $ has a
unique solution%
\begin{equation*}
u\in W_{p,\gamma }^{\left[ 2\right] }\left( 0,a;E\left( A\right) ,E\right)
\end{equation*}
for all \ $f\in L_{p}\left( 0,a;E\right) $, for $\left\vert \arg \lambda
\right\vert \leq \varphi $ with sufficiently large $\left\vert \lambda
\right\vert $ and the uniform coercive estimate holds

\begin{equation*}
\sum\limits_{i=0}^{2}\left\vert \lambda \right\vert ^{1-\frac{i}{2}%
}\left\Vert u^{\left[ i\right] }\right\Vert _{L_{p}\left( 0,a;E\right)
}+\left\Vert Au\right\Vert _{L_{p}\left( 0,a;E\right) }\leq C\left\Vert
f\right\Vert _{L_{p}\left( 0,a;E\right) }.
\end{equation*}

\textbf{Proof. }Consider the transformed problem $\left( 2.6\right) $. By
the substitution 
\begin{equation}
y=a^{1-\nu }-\left( \nu -1\right) \tau ,\text{ }\tau =\tau \left( y\right) =%
\frac{1}{\nu -1}\left( a^{1-\nu }-y\right)  \tag{2.7}
\end{equation}%
the spaces $L_{p,\tilde{\gamma}}\left( -\infty ,0;E\right) ,$ $W_{p,\tilde{%
\gamma}}^{2}\left( -\infty ,0;E\left( A\right) ,E\right) $ are mapped
isomorphically onto weighted spaces%
\begin{equation*}
L_{p,\bar{\gamma}}\left( a^{1-\nu },\infty ;E\right) ,W_{p,\bar{\gamma}%
}^{2}\left( a^{1-\nu },\infty ;E\left( A\right) ,E\right) ,
\end{equation*}%
respectively, where 
\begin{equation*}
\bar{\gamma}=\gamma \left( \tau \left( y\right) \right) =y^{\frac{\nu }{%
1-\nu }}.
\end{equation*}%
Moreover, under the substitution $\left( 2.7\right) $ the problem $\left(
2.6\right) $ is transformed into the following undegenerate problem

\begin{equation}
-u^{\left( 2\right) }\left( y\right) +Au\left( y\right) =\bar{f}\left( \tau
\right) ,  \tag{2.8}
\end{equation}%
\begin{equation*}
L_{1}u=\sum\limits_{i=0}^{m}\left[ \delta _{i}u^{\left( i\right) }\left(
a^{1-\nu }\right) +\sum\limits_{j=1}^{N}\nu _{ij}u^{\left( i\right) }\left(
y_{ij}\right) \right] =0
\end{equation*}%
considered in the weighted space $L_{p,\bar{\gamma}}\left( a^{1-\nu },\infty
;E\right) ,$ where $\bar{f}\left( y\right) =f\left( \tau \left( y\right)
\right) $ and $y_{ij}\in \left( a^{1-\nu },\infty \right) .$

By Theorem A$_{4}$ we obtain that the problem $\left( 2.8\right) $ has a
unique solution $u\in W_{p,\bar{\gamma}}^{2}\left( a^{1-\nu },\infty
;E\left( A\right) ,E\right) $ for all \ $f\in L_{p,\bar{\gamma}}\left(
a^{1-\nu },\infty ;E\right) $, $\left\vert \arg \lambda \right\vert \leq
\varphi $ with sufficiently large $\left\vert \lambda \right\vert $ and the
uniform coercive estimate holds

\begin{equation*}
\sum\limits_{i=0}^{2}\left\vert \lambda \right\vert ^{1-\frac{i}{2}%
}\left\Vert u^{\left( i\right) }\right\Vert _{L_{p,\bar{\gamma}}\left(
a^{1-\nu }\text{;}\infty ;E\right) }+\left\Vert Au\right\Vert _{_{L_{p,\bar{%
\gamma}}\left( a^{1-\nu };\infty ;E\right) }}\leq C\left\Vert f\right\Vert
_{L_{p,\bar{\gamma}}\left( a^{1-\nu }\text{;}\infty ;E\right) }.
\end{equation*}

From the above estimate we obtain the assertion.

Consider the operator $B$ generated by problem $\left( 2.4\right) $, i.e. 
\begin{equation*}
D\left( B\right) =W_{p,\gamma }^{\left[ 2\right] }\left( 0,a;E\left(
A\right) ,E,L_{1}\right) \text{, }Bu=-u^{\left[ 2\right] }+Au.
\end{equation*}

\textbf{Result 2.1. }From the Proposition 2.1 we obtain that the operator $B$
is positive in $L_{p}\left( 0,a;E\right) $ and there is positive constants $%
C_{1}$ and $C_{2}$ that 
\begin{equation}
C_{1}\left\Vert \left( B+d\right) u\right\Vert _{L_{p}\left( 0,a;E\right)
}\leq \left\Vert u\right\Vert _{W_{p,\gamma }^{\left[ 2\right] }\left(
0,a;E\left( A\right) ,E\right) }\leq C_{2}\left\Vert \left( B+d\right)
u\right\Vert _{L_{p}\left( 0,a;E\right) }  \tag{2.9}
\end{equation}%
for sufficiently large $d>0$ and $u\in D\left( B\right) .$

In a similar way as in $\left[ \text{16, Theorem 3.1}\right] $ we obtain

\textbf{Proposition 2.2}$.$ \textbf{\ }Assume the Condition 1.1 are
satisfied with $t_{1j}=\tau _{1j}$. Then, the operator $B$ is $R$-positive
in $L_{p}\left( 0,a;E\right) .$

\bigskip From $\left[ \text{15, Theorem 1}\right] $ and Remark 2.1 we obtain

\textbf{Theorem A}$_{5}$. Suppose the following conditions are satisfied:

(1) $E$ is an UMD space and\ $A$ is an $R$-positive operator in $E;$

(3) $\gamma _{k}\left( x\right) =x_{k}^{\nu _{k}}$, $\nu _{k}\in \left(
1,p\right) $, $p\in \left( 1,\infty \right) $ and $m$ is an integer, $%
\varkappa =\frac{\left\vert \alpha \right\vert }{m}\leq 1,$ $1<p<\infty $;

(4)\ \ $\Omega =\dprod\limits_{k=1}^{n}\left( 0,a_{k}\right) .$

Then, the embedding 
\begin{equation*}
D^{\alpha }W_{p,\gamma }^{\left[ m\right] }\left( \Omega ;E\left( A\right)
,E\right) \subset L_{p}\left( \Omega ;E\left( A^{1-\varkappa -\mu }\right)
\right)
\end{equation*}
is continuous. Moreover for all $h>0$ with $h\leq $ $h_{0}<\infty $ and $%
u\in W_{p,\gamma }^{m}\left( \Omega ;E\left( A\right) ,E\right) $ the
following uniform estimate holds 
\begin{equation*}
\left\Vert D^{\left[ \alpha \right] }u\right\Vert _{L_{p}\left( \Omega
;E\left( A^{1-\varkappa -\mu }\right) \right) }\leq h^{\mu }\left\Vert
u\right\Vert _{W_{p,\gamma }^{\left[ m\right] }\left( \Omega ;E\left(
A\right) ,E\right) }+h^{-\left( 1-\mu \right) }\left\Vert u\right\Vert
_{L_{p}\left( \Omega ;E\right) }.
\end{equation*}

Consider now, the following degenerate problem

\begin{equation}
\ -t^{2\nu }u^{\left( 2\right) }\left( t\right) +\left( A+\lambda \right)
u\left( t\right) =f,\text{ }t\in \left( 0,a\right) ,  \tag{2.10}
\end{equation}%
\begin{equation*}
L_{1}u=\sum\limits_{i=0}^{m}\left[ \delta _{i}u^{\left[ i\right] }\left(
a\right) +\sum\limits_{j=1}^{N}\nu _{ij}u^{\left[ i\right] }\left(
t_{ij}\right) \right] =0,
\end{equation*}%
where $u^{\left[ i\right] }=\left( t^{\nu }\frac{d}{dt}\right) ^{i}$, $\nu
>1,$ $m\in \left\{ 0,1\right\} ,$ $\delta _{i},$ $\nu _{ij},$ are complex
numbers, $t_{ij}\in \left( 0,a\right) $ and $A$ is a linear operator in $E.$

Here, 
\begin{equation*}
\gamma \left( t\right) =t^{2\nu }.
\end{equation*}

\textbf{Proposition 2.3}$.$ Assume all conditions of the Proposition 2.1 are
satisfied. Then, problem $\left( 2.10\right) $ has a unique solution $u\in
W_{p,\gamma }^{2}\left( 0,a;E\left( A\right) ,E\right) $ for all \ $f\in
L_{p}\left( 0,a;E\right) $, $\left\vert \arg \lambda \right\vert \leq
\varphi $ and sufficiently large $\left\vert \lambda \right\vert .$
Moreover, the uniform coercive estimate holds

\begin{equation}
\sum\limits_{i=0}^{2}\left\vert \lambda \right\vert ^{1-\frac{i}{2}%
}\left\Vert t^{i\nu }u^{\left( i\right) }\right\Vert _{L_{p}\left(
0,a;E\right) }+\left\Vert Au\right\Vert _{L_{p}\left( 0,a;E\right) }\leq
C\left\Vert f\right\Vert _{L_{p}\left( 0,a;E\right) }.  \tag{2.11}
\end{equation}

\textbf{Proof. \ }Since $\nu >1,$ by Theorem A$_{5}$ for all $\varepsilon >0$
there is a continuous function $C\left( \varepsilon \right) $ such that 
\begin{equation}
\left\Vert x^{\nu -1}u^{\left[ 1\right] }\right\Vert _{L_{p}\left(
0,a;E\right) }\leq \varepsilon \left\Vert u\right\Vert _{W_{p,\gamma }^{%
\left[ 2\right] }\left( 0,a;E\left( A\right) ,E\right) }+C\left( \varepsilon
\right) \left\Vert u\right\Vert _{L_{p}\left( 0,a;E\right) }  \tag{2.12}
\end{equation}%
for all $u\in W_{p,\gamma }^{\left[ 2\right] }\left( 0,a;E\left( A\right)
,E\right) $. By Result 2.1, the operator $B$ is positive in $L_{p}\left(
0,a;E\right) $. Then, in view of $\left( 2.9\right) ,$ $\left( 2.12\right) $%
, by Proposition 2.1 and resolvent properties of positive operator (see
Definition1) we have 
\begin{equation*}
\left\Vert \nu x^{\nu -1}u^{\left[ 1\right] }\right\Vert _{L_{p}\left(
0,a;E\right) }\leq \varepsilon \left\Vert \left( B+\lambda \right)
u\right\Vert _{L_{p}\left( 0,a;E\right) }+C\left( \varepsilon \right)
\left\Vert u\right\Vert _{L_{p}\left( 0,a;E\right) }\leq
\end{equation*}%
\begin{equation*}
\varepsilon \left\Vert \left( B+\lambda \right) u\right\Vert _{L_{p}\left(
0,a;E\right) }+\frac{C\left( \varepsilon \right) }{\left\vert \lambda
\right\vert }\left\Vert \left( B+\lambda \right) u\right\Vert _{L_{p}\left(
0,a;E\right) }
\end{equation*}%
for each $u\in D\left( B\right) $. From the above estimate we obtain 
\begin{equation}
\left\Vert \nu x^{\nu -1}u^{\left[ 1\right] }\right\Vert _{L_{p}\left(
0,a;E\right) }<\delta \left\Vert \left( B+\lambda \right) u\right\Vert
_{L_{p}\left( 0,a;E\right) },  \tag{2.13}
\end{equation}%
where \ $\delta =\varepsilon +\frac{C\left( \varepsilon \right) }{\left\vert
\lambda \right\vert }<1$ for sufficiently large $\left\vert \lambda
\right\vert >0.$ S\i nce $-x^{2\nu }u^{\left( 2\right) }=-u^{\left[ 2\right]
}+\nu x^{\nu -1}u^{\left[ 1\right] },$ the assertion is obtained from
Proposition 2.1 and estimate $\left( 2.13\right) .$

Consider the operator $S$ generated by problem $\left( 2.10\right) $, i.e. 
\begin{equation*}
D\left( S\right) =W_{p,\gamma }^{2}\left( 0,a;E\left( A\right)
,E,L_{k}\right) ,\text{ }Su=-x^{2\nu }u^{\left( 2\right) }+Au.
\end{equation*}

\ \textbf{Result 2.2}. Suppose all conditions of Proposition 2.1 are
satisfied. Then, the operator $S$ is $R$-positive in $L_{p}\left(
0,a;E\right) .$

\textbf{Proof}. Indeed, by Proposition 2.2 the operator $B$ is $R$-positive
in $L_{p}\left( 0,a;E\right) $. By definition of $R$ positive operators (see
Definition 3) 
\begin{equation}
\text{ }R\left\{ \lambda \left( B+\lambda \right) ^{-1},\text{ }\lambda \in
S_{\varphi }\right\} \leq M_{1}\text{.}  \tag{2.14}
\end{equation}

\textbf{\ }Then by estimates\textbf{\ }$\left( 2.13\right) $, $\left(
2.14\right) ,$ definition of $R$-bounded sets (see Definition 2) and in view
of the Kahane's contraction principle and from the product properties of the
collection of $R$-bounded operators (see e.g. $\left[ \text{4}\right] $
Lemma 3.5, Proposition 3.4) we obtain 
\begin{equation*}
\text{ }R\left\{ \lambda \left( S+\lambda \right) ^{-1},\text{ }\lambda \in
S_{\varphi }\right\} \leq M_{2}\text{.}
\end{equation*}

Consider now the leading part of the problem $\left( 2.1\right) -\left(
2.2\right) $, i.e. 
\begin{equation}
-\dsum\limits_{k=1}^{n}x_{k}^{2\alpha _{k}}\frac{\partial ^{2}u}{\partial
x_{k}^{2}}+Au+\lambda u=f\left( x\right) ,\text{ }L_{k}u=0\text{, }%
k=1,2,...,n.  \tag{2.15}
\end{equation}

\textbf{Proposition 2.4. }Assume the Condition 2.1 are satisfied. Then
problem $\left( 2.15\right) $ has a unique solution $u\in W_{p,\alpha
}^{2}\left( G;E\left( A\right) ,E\right) $ for $f\in L_{p}\left( G;E\right) $
and sufficiently large $\left\vert \lambda \right\vert $ with $\left\vert
\arg \lambda \right\vert \leq \varphi .$ Moreover, the uniform coercive
estimate holds

\begin{equation}
\sum\limits_{k=1}^{n}\sum\limits_{i=0}^{2}\left\vert \lambda \right\vert ^{1-%
\frac{i}{2}}\left\Vert x_{k}^{i\alpha }\frac{\partial ^{i}u}{\partial
x_{k}^{i}}\right\Vert _{L_{p}\left( G;E\right) }+\left\Vert Au\right\Vert
_{L_{p}\left( G;E\right) }\leq M\left\Vert f\right\Vert _{L_{p}\left(
G;E\right) }.  \tag{2.16}
\end{equation}

\textbf{Proof. }Consider first, the problem $\left( 2.15\right) $ for $n=2$
i.e

\begin{equation}
-\dsum\limits_{k=1}^{2}x_{k}^{2\alpha _{k}}\frac{\partial ^{2}u}{\partial
x_{k}^{2}}+Au+\lambda u=f\left( x_{1},x_{2}\right) ,\text{ }L_{k}u=0\text{, }%
k=1,2.  \tag{2.17}
\end{equation}%
Since%
\begin{equation*}
L_{p}\left( 0,a_{2};L_{p}\left( 0,a_{1};E\right) \right) =L_{p}\left( \left(
0,a_{1}\right) \left( 0,a_{2}\right) \times ;E\right) ,
\end{equation*}%
then the BVP $\left( 2.17\right) $ can be expressed as

\begin{equation}
-x^{2\alpha _{2}}\frac{d^{2}u}{dx_{2}^{2}}+\left( S+\lambda \right) u\left(
x_{2}\right) =f\left( x_{2}\right) \text{, }L_{2}u=0.  \tag{2.18}
\end{equation}

By virtue of $\left[ \text{1, Theorem 4.5.2}\right] $, $F=L_{p}\left(
0,a_{1};E\right) \in UMD$ provided $E\in UMD$ for $p\in \left( 1,\infty
\right) $. By Result 2.2\ the operator $S\ $is $R$-positive in $F.$ Then by
virtue of Proposition 2.3 we get that, for $f\in L_{p}\left(
0,a_{2};F\right) $ the problem $\left( 2.18\right) ,$ i.e. problem $\left(
2.17\right) $ for $\left\vert \arg \lambda \right\vert \leq \varphi $ and
sufficiently large $\left\vert \lambda \right\vert $ has a unique solution $%
u\in $ $W_{p,\alpha _{2}}^{2}\left( 0,a_{2};D\left( S\right) ,F\right) $ and
the coercive uniform estimate $\left( 2.16\right) $ holds for solution of
the problem $\left( 2.12\right) $. By continuing the above proses $n$ time,
we obtain that problem $\left( 2.15\right) $\ has a unique solution $u\in
W_{p,\alpha }^{2}\left( G;E\left( A\right) ,E\right) $ for $f\in L_{p}\left(
G;E\right) $, $\left\vert \arg \lambda \right\vert \leq \varphi $ and
sufficiently large $\left\vert \lambda \right\vert ,$ moreover, the uniform
estimate $\left( 2.16\right) $ holds.

\textbf{Proof of Theorem 2.1.} Let $Q_{0}$ denote the operator generated by
problem $\left( 2.15\right) $ i.e., 
\begin{equation*}
D\left( Q_{0}\right) =\left\{ u\in W_{p,\alpha }^{2}\left( G;E\left(
A\right) ,E\right) ,\text{ }L_{k}u=0,\text{ }k=1,2,...,n\right\} ,
\end{equation*}%
\begin{equation*}
Q_{0}u=-\dsum\limits_{k=1}^{n}x_{k}^{2\alpha _{k}}\frac{\partial ^{2}u}{%
\partial x_{k}^{2}}+Au.
\end{equation*}%
The estimate $\left( 2.16\right) $ implies that the operator $Q_{0}$ has a
bounded inverse from $L_{p}\left( G;E\right) $ to $W_{p,\alpha }^{2}\left(
G;E\left( A\right) ,E\right) $, i.e. the following estimate holds%
\begin{equation*}
\left\Vert \left( Q_{0}+\lambda \right) ^{-1}f\right\Vert _{W_{p,\alpha
}^{2}\left( G;E\left( A\right) ,E\right) }\leq C\left\Vert f\right\Vert
_{L_{p}\left( G;E\right) }
\end{equation*}%
for all $f\in L_{p}\left( G;E\right) ,$ $\lambda \in S\left( \varphi \right) 
$ with sufficiently large $\left\vert \lambda \right\vert $. Moreover, by
Theorem A$_{1}$ and in view of assumption (3), for all $\varepsilon >0$
there is a continuous function $C\left( \varepsilon \right) $ such that 
\begin{equation*}
\dsum\limits_{k=1}^{n}\left\Vert x_{k}^{\alpha _{k}}A_{k}u\right\Vert
_{L_{p}\left( G;E\right) }\leq \varepsilon \left\Vert u\right\Vert
_{W_{p,\alpha }^{2}\left( G;E\left( A\right) ,E\right) }+C\left( \varepsilon
\right) \left\Vert u\right\Vert _{L_{p}\left( G;E\right) }.
\end{equation*}%
From the above estimates we obtain that there is a positive number $\delta
<1 $ such that 
\begin{equation*}
\left\Vert Q_{1}u\right\Vert _{L_{p}\left( G;E\right) }<\delta \left\Vert
\left( Q_{0}+\lambda \right) u\right\Vert _{L_{p}\left( G;E\right) }
\end{equation*}%
for $u\in W_{p,\alpha }^{2}\left( G;E\left( A\right) ,E\right) ,$ where 
\begin{equation*}
Q_{1}u=\dsum\limits_{k=1}^{n}x_{k}^{\alpha _{k}}A_{k}\left( x\right) \frac{%
\partial u}{\partial x_{k}}.
\end{equation*}

Let $Q$ denote differential operator generated by problem $\left( 2.1\right)
-\left( 2.2\right) $ for $\lambda =0.$ It is clear that 
\begin{equation*}
\left( Q+\lambda \right) =\left[ I+Q_{1}\left( Q_{0}+\lambda \right) ^{-1}%
\right] \left( Q_{0}+\lambda \right) .
\end{equation*}%
Therefore, we obtain that the operator $\left( Q+\lambda \right) ^{-1}$ is
bounded from $L_{p}\left( G;E\right) $ to $W_{p,\alpha }^{2}\left( G;E\left(
A\right) ,E\right) $ and the estimate $\left( 2.16\right) $ is satisfied.

Let $L=L\left( L_{p}\left( G;E\right) \right) .$ We get the following result
from Theorem 2.1:

\textbf{Result 2.3. }Theorem 2.1 implies that differential operator $Q$ has
a resolvent $\left( Q+\lambda \right) ^{-1}$ for $\left\vert \arg \lambda
\right\vert \leq \varphi $ and the following estimate holds 
\begin{equation*}
\sum\limits_{k=1}^{n}\sum\limits_{i=0}^{2}\left\vert \lambda \right\vert ^{1-%
\frac{i}{2}}\left\Vert x_{k}^{i\alpha }\frac{\partial ^{i}}{\partial
x_{k}^{i}}\left( Q+\lambda \right) ^{-1}\right\Vert _{L}+\left\Vert A\left(
Q+\lambda \right) ^{-1}\right\Vert _{L}\leq M.
\end{equation*}

\begin{center}
\textbf{3. Spectral properties of singular degenerate elliptic operators }
\end{center}

In this section, the spectral properties for singular degenerate abstract
differential operators are derived. Note that, the leading part of this
operator is non-self-adjoint. Consider the differential operator $Q$
generated by BVP $\left( 2.1\right) -\left( 2.2\right) $ for $\lambda =0$.
Let%
\begin{equation*}
X=L_{p}\left( G;E\right) \text{ and }Y=W_{p\mathbf{,}\alpha }^{2}\left(
G;E\left( A\right) ,E\right) .
\end{equation*}%
The main results of this section are the following theorems:

\textbf{Theorem 3.1. }Assume the conditions of Theorem 2.1 are satisfied and 
$A^{-1}$ is compact in $E$. Then, problem $\left( 2.1\right) -\left(
2.2\right) $ is Fredholm in $L_{p}\left( G;E\right) $ for $\lambda =0.$

\textbf{Proof. }Theorem 2.1 implies that the operator $Q\mathbf{+}\lambda $
has a bounded inverse $\left( Q\mathbf{+}\lambda \right) ^{-1}$ from $%
L_{p}\left( G;E\right) $ to $W_{p\mathbf{,}\alpha }^{2}\left( G;E\left(
A\right) ,E\right) $ for\ sufficiently large $\left\vert \lambda \right\vert
,$ that is the operator $Q\mathbf{+}\lambda $ is Fredholm from $W_{p\mathbf{,%
}\alpha }^{2}\left( G;E\left( A\right) ,E\right) $ into $L_{p}\left(
G;E\right) $. Then by Theorem A$_{2}$ and in view of perturbation theory of
linear operators we obtain that the operator $Q$ is Fredholm from $W_{p%
\mathbf{,}\alpha }^{2}\left( G;E\left( A\right) ,E\right) $ into $%
L_{p}\left( G;E\right) $.

\textbf{Theorem 3.2. }Suppose all conditions of Theorem 3.1 hold, $\alpha
_{k}<2$ and 
\begin{equation*}
\text{ }s_{j}\left( I\left( E\left( A\right) ,E\right) \right) \sim j^{-%
\frac{1}{\nu }},\text{ }j=1,2,...,\infty ,\text{ }\nu >0.
\end{equation*}%
Then:

(a) 
\begin{equation}
s_{j}\left( \left( Q+\lambda \right) ^{-1}\left( L_{p}\left( G;E\right)
\right) \right) \sim j^{-\frac{1}{\nu +\varkappa }},  \tag{3.1}
\end{equation}%
where 
\begin{equation*}
\varkappa =\dsum\limits_{k=1}^{n}\frac{1}{2-\alpha _{k}}.
\end{equation*}

(b) the system of root functions of operator $Q$ is complete in $X.$

\textbf{Proof.} By virtue Theorem 4.1$,$ there exists a resolvent operator $%
\left( Q+\lambda \right) ^{-1}$ which is bounded from $X$ to $Y.$ Moreover,
by virtue of Theorem A$_{3}$ the embedding operator $I\left( Y,X\right) $ is
compact and 
\begin{equation}
s_{j}\left( I\left( Y,X\right) \right) \sim j^{-\frac{1}{\nu +\varkappa }}. 
\tag{3.2}
\end{equation}

It is clear to see that 
\begin{equation}
\left( Q+\lambda \right) ^{-1}\left( X\right) =\left( Q+\lambda \right)
^{-1}\left( X,Y\right) \times I\left( Y,X\right) .  \notag
\end{equation}%
Hence, from the relation $\left( 3.1\right) $ and Theorem A$_{3}$ we obtain $%
\left( 3.2\right) $. The Result 2.3 and the relation $\left( 3.2\right) $\
implies that operator $Q$ $+\lambda _{0}$ is positive in $X$ for
sufficiently large $\lambda _{0}$ and 
\begin{equation}
\left( Q+\lambda _{0}\right) ^{-1}\in \sigma _{q}\left( X\right) ,\text{ }%
q>\nu +\varkappa .  \tag{3.3}
\end{equation}

Then in view of the Result 2.3, the relation $\left( 3.3\right) $ and by
virtue of $\left[ \text{2, Theorem 3.4.1}\right] $ we obtain the assertion
(b).

\begin{center}
\textbf{\ 4. Singular degenerate boundary value problems for infinite
systems of equations }
\end{center}

Consider the infinite system of BVPs%
\begin{equation}
\dsum\limits_{k=1}^{n}-x_{k}^{2\alpha _{k}}\frac{\partial ^{2}u_{m}}{%
\partial x_{k}^{2}}+\dsum\limits_{k=1}^{n}\dsum\limits_{j=1}^{\infty
}x_{k}^{\alpha _{k}}a_{kmj}\left( x,y\right) \frac{\partial u_{j}}{\partial
x_{k}}+  \tag{4.1}
\end{equation}%
\begin{equation*}
\lambda u=f_{m}\left( x,y\right) ,
\end{equation*}%
\begin{equation}
\text{ }L_{1}u=0,\text{ }L_{2}u=0,  \tag{4.2}
\end{equation}%
where $\lambda $ is a complex parameter, $L_{k}$ are defined by $\left(
2.2\right) $ and $x\in G=\dprod\limits_{k=1}^{n}\left( 0,a_{k}\right) $. 
\begin{equation*}
D=\left\{ d_{m}\right\} ,\text{ }d_{m}>0,\text{ }u=\left\{ u_{m}\right\} ,%
\text{ }Du=\left\{ d_{m}u_{m}\right\} ,\text{ }m=1,2,...,
\end{equation*}

\begin{equation*}
\text{ }l_{q}\left( D\right) =\left\{ u\text{: }u\in l_{q},\right.
=\left\Vert u\right\Vert _{l_{q}\left( D\right) }=\left. \left(
\sum\limits_{m=1}^{\infty }\left\vert d_{m}u_{m}\right\vert ^{q}\right) ^{%
\frac{1}{q}}<\infty ,q\in \left( 1,\infty \right) \right\} .
\end{equation*}%
Let $O$ denote the operator in $L_{p}\left( G;l_{q}\right) $ generated by
problem $\left( 4.1\right) -\left( 4.2\right) .$ Here,%
\begin{equation*}
\alpha _{k}\left( x\right) =x_{k}^{2\alpha _{k}},\text{ }\alpha =\alpha
\left( x\right) =\left( x_{1}^{2\alpha _{1}},x_{2}^{2\alpha
_{2}},...,x_{n}^{2\alpha _{n}}\right) .
\end{equation*}%
From Theorem 2.1, we obtain

\textbf{Theorem 4.1. }Assume $\alpha _{k}\in \left( 1,p\right) $ for $p\in
\left( 1,\infty \right) $, $k=1,2,...,n\ $\ and $a_{kmj}\in L_{\infty
}\left( G\right) $. Moreover, for $0<\mu _{k}<\frac{1}{2}$\ and for all $%
x\in G$ 
\begin{equation*}
\text{ }\sup\limits_{m}\sum\limits_{j=1}^{\infty }a_{kmj}\left( x,y\right)
d_{j}^{-\left( \frac{1}{2}-\mu _{k}\right) }<M_{k}.
\end{equation*}%
Then:

(a) for all $f\left( x\right) =\left\{ f_{m}\left( x\right) \right\}
_{1}^{\infty }\in L_{p}\left( G;l_{q}\right) ,$ $p,q\in \left( 1,\infty
\right) $, $\left\vert \arg \lambda \right\vert \leq \varphi $, $0\leq
\varphi <\pi $ and for sufficiently large $\left\vert \lambda \right\vert $
problem $\left( 4.1\right) -\left( 4.2\right) $ has a unique solution $%
u=\left\{ u_{m}\left( x\right) \right\} _{1}^{\infty }$ that belongs to $%
W_{p,\alpha }^{2}\left( G,l_{q}\left( D\right) ,l_{q}\right) $ and 
\begin{equation*}
\dsum\limits_{k=1}^{n}\sum\limits_{i=0}^{2}\left\vert \lambda \right\vert
^{1-\frac{i}{2}}\left\Vert x_{k}^{i\alpha }\frac{\partial ^{i}u}{\partial
x_{k}^{i}}\right\Vert _{L_{p}\left( G;l_{q}\right) }+\left\Vert
Du\right\Vert _{L_{p}\left( G;l_{q}\right) }\leq M\left\Vert f\right\Vert
_{L_{p}\left( G;l_{q}\right) };
\end{equation*}

(b) the operator $O$ is Fredholm in $L_{p}\left( G;l_{q}\right) ;$

(c) the system of root functions of operator $O$ is complete in $L_{p}\left(
G;l_{q}\right) .$

\textbf{Proof. } Let $E=l_{q},$ $A$ and $A_{\alpha }\left( x\right) $ be
infinite matrices, such that

\begin{equation*}
A=\left[ d_{m}\delta _{mj}\right] ,\text{ }A_{k}\left( x\right) =\left[
d_{kjm}\left( x\right) \right] ,\text{ }m\text{, }j=1,2,...\infty .
\end{equation*}%
By $\left[ 4\right] ,$ $l_{q}$\ is the UMD space. It is clear to see that
the operator $A$ is $R$-positive in $l_{q}$. The problem $\left( 4.1\right) $
can be rewritten in the form of $\left( 2.1\right) -\left( 2.2\right) $.
From Theorem 2.1 we obtain that problem $\left( 4.1\right) -\left(
4.2\right) $ has a unique solution $u\in $ $W_{p,\alpha }^{2}\left(
G;l_{q}\left( D\right) ,l_{q}\right) $ for all $f\in L_{p}\left(
G;l_{q}\right) $ and 
\begin{equation*}
\sum\limits_{k=1}^{n}\sum\limits_{i=0}^{2}\left\vert \lambda \right\vert ^{1-%
\frac{i}{2}}\left\Vert x_{k}^{i\alpha }\frac{\partial ^{i}u}{\partial
x_{k}^{i}}\right\Vert _{L_{p}\left( G;l_{q}\right) }+\left\Vert
Du\right\Vert _{L_{p}\left( G;l_{q}\right) }\leq M\left\Vert f\right\Vert
_{L_{p}\left( G;l_{q}\right) }.
\end{equation*}

From the above estimate we obtain the assertion (a). The assertions (b) and
(c) are obtained from Theorems 3.1 and 3.2, respectively.

\begin{center}
\textbf{5. Wentzell-Robin type BVP for degenerate elliptic equation }
\end{center}

Consider the problem%
\begin{equation}
\dsum\limits_{k=1}^{n}x_{k}^{2\alpha _{k}}\frac{\partial ^{2}u}{\partial
x_{k}^{2}}+a\frac{\partial ^{2}u}{\partial y^{2}}+b\frac{\partial u}{%
\partial y}+\lambda u=f\left( x,y\right) ,\text{ }x\in G,\text{ }y\in \left(
0,1\right) ,  \tag{5.1}
\end{equation}

\begin{equation}
L_{k}u=0\text{, for a.e. }y\in \left( 0,1\right) ,  \tag{5.2}
\end{equation}

\begin{equation}
\text{ }a\left( j\right) u_{yy}\left( x,j\right) +b\left( j\right)
u_{y}\left( x,j\right) =0=0\text{, }j=0,1,\text{ for a.e. }x\in G,  \tag{5.3}
\end{equation}%
where $a=a\left( y\right) ,$ $b=b\left( y\right) $ are real-valued functions
on $\left( 0,1\right) $, $\lambda $ is acomplex parameter and $L_{k}$ are
boundary condition in $x$ defined by $\left( 2.2\right) .$ For $\Omega
=G\times \left( 0,1\right) $, $\mathbf{p=}\left( p,2\right) $ and $L_{%
\mathbf{p}}\left( \Omega \right) $ will denote the space of all $\mathbf{p}$%
-summable scalar-valued\ functions with mixed norm. Analogously, $W_{\mathbf{%
p,}\alpha }^{2}\left( \Omega \right) $ denotes the Sobolev space with
corresponding mixed norm, i.e., $W_{\mathbf{p,}\alpha }^{2}\left( \Omega
\right) $ denotes the space of all functions $u\in L_{\mathbf{p}}\left(
\Omega \right) $ possessing the derivatives $x_{k}^{2\alpha _{k}}\frac{%
\partial ^{2}u}{\partial x_{k}^{2}}\in L_{\mathbf{p}}\left( \Omega \right) $
with the norm 
\begin{equation*}
\ \left\Vert u\right\Vert _{W_{\mathbf{p,}\alpha }^{2}\left( \Omega \right)
}=\left\Vert u\right\Vert _{L_{\mathbf{p}}\left( \Omega \right)
}+\dsum\limits_{k=1}^{n}\left\Vert x^{2\alpha _{k}}\frac{\partial ^{2}u}{%
\partial x_{k}^{2}}\right\Vert _{L_{\mathbf{p}}\left( \Omega \right) }.
\end{equation*}

\bigskip \textbf{Condition 5.1 }Assume;

(1) $\alpha _{k}\in \left( 1,p\right) $ for $p\in \left( 1,\infty \right) $
and $k=1,2,...,n;$

(2)\ $a$ is positive, $b$ is a real-valued functions on $\left( 0,1\right) ;$

(3) $a\left( .\right) \in C\left[ 0,1\right] $ and%
\begin{equation*}
\exp \left( -\dint\limits_{\frac{1}{2}}^{x}b\left( t\right) a^{-1}\left(
t\right) dt\right) \in L_{1}\left( 0,1\right) .
\end{equation*}

\bigskip Let $H$ denote the elliptic operator in $L_{\mathbf{p}}\left(
\Omega \right) $ generated by problem $\left( 5.1\right) -\left( 5.3\right)
. $\ In this section, we present the following result:

\bigskip \bigskip \textbf{Theorem 5.1. }Suppose the Condition 5.1 hold. Then:

(a) for $f\in L_{\mathbf{p}}\left( \Omega \right) $ problem $\left(
5.1\right) -\left( 5.3\right) $ for $\lambda \in S\left( \varphi \right) $
and sufficiently large $\left\vert \lambda \right\vert >0$\ has a unique
solution $u$\ belonging to $W_{\mathbf{p,}\alpha }^{2}\left( \Omega \right) $
and the following coercive uniform estimate holds 
\begin{equation*}
\sum\limits_{k=1}^{n}\sum\limits_{i=0}^{2}\left\vert \lambda \right\vert ^{1-%
\frac{i}{2}}\left\Vert x_{k}^{i\alpha _{k}}\frac{\partial ^{i}u}{\partial
x_{k}^{i}}\right\Vert _{L_{\mathbf{p}}\left( \Omega \right) }\leq
C\left\Vert f\right\Vert _{L_{\mathbf{p}}\left( \Omega \right) };
\end{equation*}

(b) the problem $\left( 5.1\right) -\left( 5.3\right) $ is Fredholm in $L_{%
\mathbf{p}}\left( \Omega \right) $ for $\lambda =0;$

(c) the system of root functions of operator $H$ is complete in $L_{p}\left(
G;l_{q}\right) $

\ \textbf{Proof.} Let $E=L_{2}\left( 0,1\right) $. It is known $\left[ 10%
\right] $\ that $L_{2}\left( 0,1\right) $ is an $UMD$ space. Consider the
operator $A$ defined by 
\begin{equation*}
D\left( A\right) =W_{2}^{2}\left( 0,1;A\left( j\right) u\left( j\right)
=0\right) ,\text{ }j=0,1,\text{ }Au=a\frac{\partial ^{2}u}{\partial y^{2}}+b%
\frac{\partial u}{\partial y}.
\end{equation*}

Therefore, the problem $\left( 5.1\right) -\left( 5.2\right) $ can be
rewritten in the form of $\left( 2.1\right) -\left( 2.2\right) $, where $%
u\left( x\right) =u\left( x,.\right) ,$ $f\left( x\right) =f\left(
x,.\right) $\ are functions with values in $E=L_{2}\left( 0,1\right) .$ By
virtue of $\left[ \text{9, 10}\right] $ the operator $A$ generates analytic
semigroup in $L_{2}\left( 0,b\right) $. Then in view of Hill-Yosida theorem
(see e.g. $\left[ \text{20, \S\ 1.13}\right] $) this operator is positive in 
$L_{2}\left( 0,b\right) .$ Since all uniform bounded set in Hilbert space is 
$R$-bounded (see $\left[ 4\right] $ ), i.e. we get that the operator $A$ is $%
R$-positive in $L_{2}\left( 0,b\right) .$ Then from Theorem 2.1 we obtain
the assertion (a). Since the embedding $W_{2}^{2}\left( 0,1\right) \subset
L_{2}\left( 0,1\right) $ is compact, the assertions (b) and (c) are obtained
from Theorems 3.1 and 3.2, respectively.

From Theorem 5.1 we obtain:

\textbf{Result 5.1. }Theorem 5.1 implies that operator $H$ has a resolvent $%
\left( H+\lambda \right) ^{-1}$ for $\left\vert \arg \lambda \right\vert
\leq \varphi $ and the following sharp coercive resolvent estimate holds 
\begin{equation*}
\sum\limits_{k=1}^{n}\sum\limits_{i=0}^{2}\left\vert \lambda \right\vert ^{1-%
\frac{i}{2}}\left\Vert x_{k}^{i\alpha }\frac{\partial ^{i}}{\partial
x_{k}^{i}}\left( H+\lambda \right) ^{-1}\right\Vert _{L\left( L_{\mathbf{p}%
}\left( \Omega \right) \right) }\leq M.
\end{equation*}

\bigskip \textbf{Acknowledgements}

The author is thanking the library manager of Okan University Kenan Oztop
for his help in finding the necessary articles and books in my research area.

\bigskip \textbf{References}

\begin{quote}
\ \ \ \ \ \ \ \ \ \ \ \ \ \ \ \ \ \ \ \ \ \ \ \ \ \ \ \ \ \ \ \ \ \ \ \ \ \
\ \ \ \ \ \ \ \ \ \ \ \ \ \ \ \ \ \ \ \ \ \ \ \ \ \ \ \ \ \ \ \ \ \ \ \ \ \
\ \ \ \ \ \ \ \ \ \ \ \ \ \ \ \ \ 
\end{quote}

\begin{enumerate}
\item Amann H., Linear and quasi-linear equations,1, Birkhauser, Basel 1995.

\item Agranovicn M. S., Spectral Boundary Value Problems in Lipschitz
Domains or Strongly Elliptic Systems in Banach Spaces $H_{p}^{\sigma }$ and $%
B_{p}^{\sigma }$, Func. Anal. Appl., 42 (4), (2008), 249-267.

\item Ashyralyev A., Cuevas C., and Piskarev S., On well-posedness of
difference schemes for abstract elliptic problems in spaces, Numer. Func.
Anal.Opt., v. 29, No. 1-2, Jan. 2008, 43-65.

\item Burkholder D. L., A geometrical conditions that implies the existence
certain singular integral of Banach space-valued functions, Proc. conf.
Harmonic analysis in honor of Antonu Zigmund, Chicago, 1981,Wads Worth,
Belmont, (1983), 270-286.

\item Campiti, M., Metafune, G., and Pallara, D., Degenerate self adjoint
equations on the unit \i nterval, Semigroup Forum 57 (1998), 1 -- 36.

\item Dore C. and Yakubov S., Semigroup estimates and non coercive boundary
value problems, Semigroup Forum 60 (2000), 93-121.

\item Denk R., Hieber M., Pr\"{u}ss J., $R$-boundedness, Fourier multipliers
and problems of elliptic and parabolic type, Mem. Amer. Math. Soc. 166
(2003), n.788.

\item Favini A., Shakhmurov V., Yakubov Y., Regular boundary value problems
for complete second order elliptic differential-operator equations in UMD
Banach spaces, Semigroup Forum, v. 79 (1), 2009.

\item Favini A., Goldstein G. R., Goldstein J. A. and Romanelli S.,
Degenerate second order differential operators generating analytic
semigroups in $L_{p}$ and $W^{1,p}$, Math. Nachr. 238 (2002), 78 -- 102.

\item Keyantuo V., Lizama, C., Maximal regularity for a class of
integro-differential equations with infinite delay in Banach spaces, Studia
Math. 168 (2005), 25-50.

\item V. Keyantuo, M. Warma, The wave equation with Wentzell--Robin boundary
conditions on Lp-spaces, J. Differential Equations 229 (2006) 680--697.

\item Lunardi A., Analytic semigroups and optimal regularity in parabolic
problems, Birkhauser, 2003.

\item Lions J. L and Peetre J., Sur une classe d'espaces d'interpolation,
Inst. Hautes Etudes Sci. Publ. Math., 19(1964), 5-68.

\item Shakhmurov V. B., Estimates of approximation numbers for embedding
operators and applications, Acta Math. Sin., 28( 9 )(2012), 1883-1896.

\item Shakhmurov V. B., Abstract capacity of regions and compact embedding
with applications, Acta Math. Sci., (31)1, 2011, 49-67.

\item Shakhmurov V. B., Degenerate differential operators with parameters,
Abstr. Appl. Anal., 2007, v. 2006, 1-27.

\item Shakhmurov V. B., Singular degenerate differential operators and
applications, Math. Slovaca 65(2015)(6), 1-24.

\item Shahmurov R., On strong solutions of a Robin problem modeling heat
conduction in materials with corroded boundary, Nonlinear Anal. Real World
Appl., 2011,13(1), 441-451.

\item Shahmurov R., Solution of the Dirichlet and Neumann problems for a
modified Helmholtz equation in Besov spaces on an annuals, J. Differential
Equations, 2010, 249(3), 526-550.

\item Triebel H., Interpolation theory, Function spaces, Differential
operators, North-Holland, Amsterdam, 1978.

\item Weis L, Operator-valued Fourier multiplier theorems and maximal $L_{p}$
regularity, Math. Ann. 319, (2001), 735-758.\ 

\item Yakubov S., and Yakubov Ya., Differential-operator Equations. Ordinary
and Partial \ Differential Equations, Chapman and Hall /CRC, Boca Raton,
2000.
\end{enumerate}

\bigskip

\end{document}